\documentclass{amsart}
\usepackage{graphicx}
\usepackage{amssymb}
\usepackage{epstopdf}
\usepackage{nicefrac}
\usepackage{enumerate}
\usepackage{hyperref}
\usepackage{float}
\usepackage{color}
\usepackage{pst-all}
\usepackage{tabularx}

\newtheorem{theorem}{THEOREM}[section]
\newtheorem{thm}{THEOREM}[section]

\newtheorem{cor}[thm]{COROLLARY}

\newtheorem{defn}[thm]{DEFINITION}

\newtheorem{lemma}[thm]{LEMMA}

\newtheorem{prop}[thm]{PROPOSITION}

\newcommand{\ds}{\displaystyle}

\newcommand{\cG}{{\mathcal G}}
\newcommand{\cH}{{\mathcal H}}

\newcommand{\cO}{{\mathcal O}}

\newcommand{\cU}{{\mathcal U}}

\newcommand{\diam}{{\rm diam}} 
 
\newcommand{\dX}{d_{\fX}} 
\newcommand{\e}{{\varepsilon}}

\newcommand{\fX}{{\mathfrak{X}}}

\newcommand{\G}{\Gamma}

\newcommand{\Homeo}{\operatorname{Homeo}}

\newcommand{\mZ}{{\mathbb Z}}

\begin{document}

\title{Essential holonomy of   Cantor actions}

\author{Steven Hurder}
\address{Steven Hurder, Department of Mathematics, University of Illinois at Chicago, 322 SEO (m/c 249), 851 S. Morgan Street, Chicago, IL 60607-7045}
\email{hurder@uic.edu}

\author{Olga Lukina}
\address{Olga Lukina, Mathematical Institute, 
  Leiden University,
P.O. Box 9512,
2300 RA Leiden,
The Netherlands}
\email{o.lukina@math.leidenuniv.nl}

\thanks{Version date: January 20, 2023}



\begin{abstract}
A group action has essential holonomy if the set of points with non-trivial holonomy has positive measure. If such an action is topologically free, then having essential holonomy is equivalent to the action not being essentially free, which means that the set of points with non-trivial stabilizer has positive measure. In this paper, we investigate the relation between the property of having essential holonomy and structure of the acting group for minimal equicontinuous actions on Cantor sets. We show that if such a group action is locally quasi-analytic and has essential holonomy, then every commutator subgroup in the group lower central series has elements with positive measure set of points with non-trivial holonomy. In particular, this gives a new proof that a minimal equicontinuous Cantor action by a nilpotent group has   no essential holonomy. We also show that the property of having essential holonomy is preserved under return equivalence and continuous orbit equivalence of minimal equicontinuous Cantor actions.
Finally, we give examples to show that the assumption on the action that it is locally quasi-analytic is necessary.
\end{abstract}

\maketitle
 
\section{Introduction}\label{sec-intro}

We say that $(\fX, \G, \Phi)$ is a \emph{Cantor action} if $\G$ is a countable group, $\fX$ is a   Cantor   space, and $\Phi \colon \G \times \fX \to \fX$ is an  action by homeomorphisms. We assume throughout this paper  that, moreover, Cantor actions are minimal and equicontinuous. This implies that a Cantor action $(\fX, \G, \Phi)$  has  a unique invariant probability measure   $\mu$.   We recall   further basic definitions and   constructions for   Cantor actions in    Section \ref{sec-basics}, as used in the formulation of our results below.

The dynamical properties  of   Cantor actions can have surprisingly subtle aspects, especially for the case where $\G$ is non-abelian, as revealed by the many   examples in the literature. One approach to classifying Cantor actions is by their dynamical properties, and to also consider the relation between the algebraic properties of the acting group $\G$ and their dynamics.
The main result of this work  makes a new contribution to this classification scheme, as it relates the algebraic properties of $\G$ and the dynamics of the action, through the study of the property that an action has   \emph{essential holonomy}; see Definition~\ref{def-essentialholonomy} below.  
 
First recall some standard notions concerning the fixed-point sets for a Cantor action $(\fX, \G, \Phi)$. We use the notation $g \cdot x = \Phi(g)(x)$, for $g \in \G$ and $x \in \fX$.
The set $\fX_g = \{ x \in \fX \mid g \cdot x = x\}$ consists of fixed points for $g$,  and the \emph{stabilizer} of a point $x \in \fX$ is the subgroup $\G_x = \{ g \in \G \mid g \cdot x = x\}$. Let 
\begin{equation}\label{eq-fpset}
\fX_\G = \bigcup_{e \ne g  \in \G}  \ \fX_{g} = \{x \in \fX \mid \G_x \ne \{e\}\}  
\end{equation}
be the   set of all points fixed by some non-identity element $g \in \G$.
 \begin{defn}
 A minimal Cantor action $(\fX, \G, \Phi)$ with invariant probability measure $\mu$ is:
 \begin{enumerate}
\item \emph{free} if $\fX_{\G}$ is empty.
\item \emph{topologically free} if $\fX_\G$ is a meager set in $\fX$.
\item \emph{essentially free}  if $\mu(\fX_\G) = 0$.
\end{enumerate}
 \end{defn}

Recall   that a Cantor action  $(\fX, \G, \Phi)$ is \emph{effective} if the only element of $\G$ that acts as the identity on $\fX$ is the identity of $\G$. 
It is elementary to show that if $\G$ is abelian, then every effective  Cantor action of $\G$ must be free. 
The topologically free Cantor actions have an important role in the study of the   $C^*$-algebras associated to  the actions, as studied for example in \cite{AS1994,BoyleTomiyama1998,Kennedy2020,LeBMB2018,Renault2008}.

  Kambites, Silva and Steinberg   showed in \cite[Theorem~4.3]{KSS2006} that the action of a group generated by finite automata on a rooted tree is topologically free if and only if it is essentially free.   Joseph   proved in \cite[Corollary 2.4]{Joseph2021} that if $\G$ has countably many subgroups, then a topologically free Cantor action of $\G$ is  essentially free.

Bergeron and Gaboriau \cite{BergeronGaboriau2004} showed that if $\G$ is a non-amenable group which is a free product of two residually finite groups, then $\G$ admits a   Cantor action which is topologically free and not essentially free. Ab\'ert and Elek \cite{AE2007} proved a similar result for finitely generated non-abelian free groups $\Gamma$.  Joseph \cite{Joseph2021} proved that any non-amenable surface group admits a continuum of pairwise non-conjugate and non-measurably isomorphic   Cantor actions which are topologically free and not essentially free.

    Examples of effective   Cantor actions which are not topologically free include some actions of branch and weakly branch groups on the boundaries of rooted trees  \cite{Grigorchuk2011,Nekrashevych2005,SilvaSteinberg2005}, some actions of nilpotent groups \cite{HL2021a}, and actions of topological full groups  \cite{CortezMedynets2016,Grigorchuk2011}.

The  work by Gr\"{o}ger and Lukina \cite{GL2019} introduced a refinement of the   notion of essentially free actions,  called  the \emph{essential holonomy} property.  
The idea is to consider   the dynamics of the action in small neighborhoods of fixed points. In place  of the set of points with trivial (resp. non-trivial) stabilizers, one considers  the set of points with trivial (resp. non-trivial) holonomy. The analog of an essentially free Cantor action is an action which has \emph{no essential holonomy}. 

    \begin{defn} \label{def-holonomy}
Let $(\fX, \G, \Phi)$  be a Cantor action. 
Say that $x \in \fX$ is   a \emph{point of non-trivial holonomy for}   $g \in \G$ if  $\Phi(g)(x) = x$, and for each open set $U \subset \fX$ with $x \in U$, there exists $y \in U$ such that $\Phi(g)(y) \ne y$. 
\end{defn}

We say that a fixed point $x \in \fX$  is a \emph{point of trivial holonomy}  for $g \in \G$, if $x$ is fixed by $\Phi(g)$, and $x$ has an open neighborhood $U_{x,g}$ where every point is fixed by  $\Phi(g)$.  We say that $x \in \fX$ is a \emph{point of trivial holonomy}, if $x$ is a point of trivial holonomy for all $g \in \G$ with $g \cdot x = x$. We say that $x \in \fX$ is a \emph{point of non-trivial holonomy}, if $x$ is a point of non-trivial holonomy for some $g \in \G$ with $g \cdot x = x$.

  Let $\fX_g^{hol}$ denote the (possibly empty)  subset of points  of non-trivial holonomy in $\fX_g$, and set  
  \begin{equation}\label{eq-hfpset}
 \fX_{\G}^{hol} = \bigcup_{g \in \G}  \ \fX_{g}^{hol} \subset \bigcup_{g \in \G}  \ \fX_{g} \subset \fX_{\G} \ .
\end{equation}
The set $\fX_\G^{hol}$ is invariant under the action of $\G$, thus $\fX_\G^{hol}$ has either $\mu$-measure 0 or 1.
The following  concept  was formulated in \cite{GL2019}:
\begin{defn}\label{def-essentialholonomy}
A   Cantor group action $(\fX,\G,\Phi)$ has \emph{essential holonomy} if the set $\fX_\G^{hol}$ of points with non-trivial holonomy has positive measure. Otherwise, it has  \emph{no essential holonomy}.
\end{defn}
    Vorobets \cite{Vorobets2012} showed that the   standard action of the Grigorchuk group on the boundary of a binary rooted tree   has only a countable set of points with non-trivial holonomy, which implies that it has no essential holonomy. 
     Gr\"oger and Lukina  proved in  \cite{GL2019} that the action of a group generated by finite automata on a rooted tree has no essential holonomy; their proof  does not require the Cantor action to be topologically free. 
  They also    gave a criterion for when   a group action on a rooted tree   has  no essential holonomy for its boundary Cantor action.


The decomposition \eqref{eq-hfpset} of $\fX_{\G}^{hol}$ as a countable union of fixed-point sets   implies that an action has essential holonomy if and only if, for at least one $g \in \G$, the set $\fX_g^{hol}$ has positive $\mu$-measure. In particular, this implies that whether or not a Cantor action has essential holonomy  is a local property. 
This remark is the basis for the following   invariance result, with details and proof   in Section~\ref{sec-invariance}.

\begin{prop} \label{prop-inv}
\begin{enumerate}
\item Suppose that Cantor actions  $(\fX_i, \G_i, \Phi_i)$, for $i =1,2$, are continuous orbit equivalent. Then  either both have  essential holonomy, or both have no essential holonomy.
\item Suppose that Cantor actions  $(\fX_i, \G_i, \Phi_i)$, for $i =1,2$, are return equivalent. Then  either both have  essential holonomy, or both have no essential holonomy.
\end{enumerate}
\end{prop}
 
For a Cantor action $(\fX, \G, \Phi)$, the assumption that $\mu$ is a probability measure, and that the action is minimal, implies that $\mu$ is a continuous measure; that is, the measure of any point $x \in \fX$ is zero. Thus,  if a Cantor action has essential holonomy, then the intersection of $\fX_g^{hol}$ with the support of $\mu$ must be an uncountable set.  This remark is the basis for the result by Joseph in \cite[Corollary 2.4]{Joseph2021} that if $\G$ has the countable subgroup property, then a topologically free Cantor action by $\G$ is essentially free, and hence has no essential holonomy.

Our main result below relates the dynamical properties of a Cantor action of a finitely-generated group $\G$, with  the lower central series 
$\G = \gamma_{1}(\G) \supset \gamma_{2}(\G) \supset \cdots \supset \gamma_{n}(\G) \supset \cdots$. 
\begin{thm}\label{thm-main0}
Let   $(\fX, \G, \Phi)$ be a locally quasi-analytic   Cantor action. If the action has essential holonomy, then for every $n \geq 1$ there exists  $\phi_n \in \gamma_{n}(\G)$ such that the action of $\phi_n$ has a positive measure set of points with non-trivial holonomy.  
\end{thm}
The proof of Theorem~\ref{thm-main0} is given  in Section~\ref{sec-proof}. We observe the following corollary for $n=2$:
\begin{cor}\label{cor-commutator}
Let   $(\fX, \G, \Phi)$ be a locally quasi-analytic   Cantor action. Let $[\G,\G] \subset \G$ be the commutator subgroup of $\G$.  If the action has essential holonomy, then   there exists  $\phi_2 \in [\G,\G]$ such that the action of $\phi_2$ has a positive measure set of points with non-trivial holonomy.  
\end{cor}

It is tempting to apply Corollary~\ref{cor-commutator} inductively, so that the conclusion of Theorem~\ref{thm-main0} applies to the derived series of $\G$. This argument does not go through, though, because while the action of the commutator subgroup $[\G,\G]$ on $\fX$ is again equicontinuous and locally quasi-analytic, it need not be minimal, and the minimality assumption on the action  is    critical for the proof of Theorem~\ref{thm-main0}.

The family of examples constructed in the proof of Theorem \ref{thm-examples} show  that the assumption that the Cantor action is locally quasi-analytic is essential for the conclusions of Theorem \ref{thm-main0} and Corollary \ref{cor-commutator}. The actions of groups $\G$ and the commutator subgroups $[\G,\G]$ in Theorem \ref{thm-examples} are minimal and not locally quasi-analytic, and the action of the group $\G$ has essential holonomy, while the action of $[\G,\G]$ has no essential holonomy.

The locally quasi-analytic property is a localized form of the topologically free property, as discussed in Section~\ref{sec-basics}. If $\G$ is a Noetherian group, then every Cantor action by $\G$ must be locally quasi-analytic by \cite[Theorem~1.6]{HL2018b}.
For a nilpotent group $\G$, the lower central series terminates, and  $\G$ is Noetherian,  so as a consequence  we obtain a novel proof of a result of Joseph which he proved for topologically free Cantor actions in \cite{Joseph2021}:
\begin{cor}\label{cor-nilp}
Let $(\fX, \G, \Phi)$  be a   Cantor action. If $\G$ is a finitely-generated nilpotent group, then $(\fX, \G, \Phi)$ has no essential holonomy. 
\end{cor}

Up until recently, it was an open problem to show that   if $\G$ is a finitely-generated amenable group, then every minimal equicontinuous Cantor action of $\G$ has no essential holonomy.  In the paper \cite{Joseph2023}, Joseph constructs an action of the wreath product of two finitely generated amenable groups which is topologically free and not essentially free. Actions in Theorem \ref{thm-examples} in our paper are not locally quasi-analytic (and therefore not topologically free) actions of infinitely generated amenable groups which have essential holonomy. It would be interesting to obtain a criterion for when an amenable group admits an action with essential holonomy. 

  Theorem~\ref{thm-main0}  suggests  that   the property of having   essential holonomy  is an interesting    invariant  of  Cantor actions, intrinsically related to the structure of the acting group, to be further explored.  
 

{\bf Acknowledgements:} The paper started in response to the question by Eduardo Scarparo whether an action of an amenable group must have no essentially holonomy. The authors also thank Eduardo for bringing to their attention the work by Matthieu Joseph \cite{Joseph2021}.

  
 \section{Cantor actions}\label{sec-basics}

We recall some   basic  properties of     Cantor actions. More details can be found in   \cite{Auslander1988,CortezPetite2008,CortezMedynets2016,HL2018b,HL2019a,LavrenyukNekrashevych2002}.

 \subsection{Basic notions}\label{subsec-basic}

Let  $(\fX,\G,\Phi)$   denote a topological  action  $\Phi \colon \G \times \fX \to \fX$. We   write $g\cdot x$ for $\Phi(g)(x)$ when appropriate.
The orbit of  $x \in \fX$ is the subset $\cO(x) = \{g \cdot x \mid g \in \G\}$. 
The action is \emph{minimal} if  for all $x \in \fX$, its   orbit $\cO(x)$ is dense in $\fX$.
The action is said to be \emph{effective}, or \emph{faithful},  if   the homomorphism $\Phi \colon \G \to \Homeo(\fX)$  is injective.
    
 An action  $(\fX,\G,\Phi)$ is \emph{equicontinuous} with respect to a metric $\dX$ on $\fX$, if for all $\e >0$ there exists $\delta > 0$, such that for all $x , y \in \fX$ and $g \in \G$ we have  that 
 $\ds  \dX(x,y) < \delta$ implies   $\dX(g \cdot x, g \cdot y) < \e$.
The property of being equicontinuous    is independent of the choice of the metric   on $\fX$ which is   compatible with the topology of $\fX$.

Recall that we assume    $(\fX,\G,\Phi)$ is  a minimal   equicontinuous action.  
We say that $U \subset \fX$  is \emph{adapted} to the action  if $U$ is a   {non-empty clopen} subset, and for any $g \in \G$, 
if $\Phi(g)(U) \cap U \ne \emptyset$ then  $\Phi(g)(U) = U$.   
Recall a basic property of equicontinuous Cantor    actions (see   \cite[Section~3]{HL2018b}).
\begin{prop}\label{prop-adpatedchain}
Let  $(\fX,\G,\Phi)$   be a  Cantor    action, and let $\dX$ be an invariant metric on $\fX$ compatible with the topology on $\fX$. Given $x \in \fX$ and $\e > 0$,  there exists an adapted clopen set $U \subset \fX$ with $x \in U$ and $\diam(U) < \e$.   
\end{prop}
\begin{cor}\label{cor-subbasis}
The adapted clopen sets form a subbasis for the topology of $\fX$.
\end{cor}

For an adapted set $U$,   the set of ``return times'' to $U$, 
 \begin{equation}\label{eq-adapted}
\G_U = \left\{g \in \G \mid g \cdot U  \cap U \ne \emptyset  \right\}  
\end{equation}
is a subgroup of   $\G$, called the \emph{stabilizer} of $U$.      
  Then for $g, g' \in \G$ with $g \cdot U \cap g' \cdot U \ne \emptyset$ we have $g^{-1} \, g' \cdot U = U$, hence $g^{-1} \, g' \in \G_U$. Thus,  the  translates $\{ g \cdot U \mid g \in \G\}$ form a finite clopen partition of $\fX$, and are in 1-1 correspondence with the quotient space $X_U = \G/\G_U$. Then $\G$ acts by permutations of the finite set $X_U$ and so the stabilizer group $\G_U \subset \G$ has finite index.  Note that this implies that if $V \subset U$ is a proper inclusion of adapted sets, then the inclusion $\G_V \subset \G_U$ is also proper.
 
 Let $U$ be an adapted set $U$ for the action  $(\fX,\G,\Phi)$, then the action of $\G_U$ restricts to an action on $U$, so we have a homomorphism $\Phi_U \colon \G_U \to \Homeo(U)$. Let $\cH_{U}$ denote the image of this action. Note that the map $\Phi_U \colon \G_U \to \cH_{U}$ is injective if the action is topologically free.

\subsection{Group chains}\label{sec-gchains}

Given a basepoint $x$, by iterating the process in Proposition \ref{prop-adpatedchain} one can always construct the following:
 
\begin{defn}\label{def-adaptednbhds}
Let  $(\fX,\G,\Phi)$   be a Cantor  action.
A properly descending chain of clopen sets $\cU = \{U_{\ell} \subset \fX  \mid \ell \geq 0\}$ is said to be an \emph{adapted neighborhood basis} at $x \in \fX$ for the action $\Phi$  if each $U_{\ell}$ is adapted to the action $\Phi$, 
    $U_{\ell +1} \subset U_{\ell}$ for all $ \ell \geq 1$,  and     $\cap  \ U_{\ell} = \{x\}$.
\end{defn}

Let    $\G_{\ell} = \G_{U_{\ell}}$ denote the stabilizer group of $U_{\ell}$ given by \eqref{eq-adapted}.
Then  we obtain a descending chain of finite index subgroups 
 $ \ds  \cG^x_{\cU} : \G = \G_0 \supset \G_1 \supset \G_2 \supset \cdots $.  Note that each $\G_{\ell}$ has finite index in $\G$, and is not assumed to be a normal subgroup.  Also note that while the intersection of the chain $\cU$ is a single point $\{x\}$, the intersection of the stabilizer groups   in  $\cG^x_{\cU}$ need not be the trivial group.
 
Next, set $X_{\ell} = \G/\G_{\ell}$ and note that  $\G$ acts transitively on the left on   $X_{\ell}$.    
The inclusion $\G_{\ell +1} \subset \G_{\ell}$ induces a natural $\G$-invariant quotient map $p_{\ell +1} \colon X_{\ell +1} \to X_{\ell}$.
 Introduce the inverse limit 
 \begin{equation} \label{eq-invlimspace}
X ~ = ~  \varprojlim \ \{p_{\ell +1} \colon X_{\ell +1} \to X_{\ell} \mid \ell > 0\} = \{(x_\ell)  = (x_0,x_1,\ldots) \mid p_{\ell+1}(x_{\ell+1}) = x_\ell\}
\end{equation}
which is a Cantor space with the Tychonoff topology. Thus elements of $X$ are infinite sequences with entries in $X_\ell$, $\ell \geq 0$. The actions of $\G$ on the factors $X_{\ell}$ induce    a minimal  equicontinuous action, denoted by  $\Phi_x \colon \G \times X \to X$, which reads
  \begin{equation}\label{eq-actioninvlim}(g, (x_\ell)) \mapsto g \cdot (x_\ell) = (g \cdot x_\ell) = (g \cdot x_0, g \cdot x_1,\ldots).\end{equation}
 For each $\ell \geq 0$, we have the ``partition coding map'' $\Theta_{\ell} \colon \fX \to X_{\ell}$ which is $\G$-equivariant.  The maps $\{\Theta_{\ell}\}$ are compatible with the   map on quotients in \eqref{eq-invlimspace}, and so define a  limit map $\Theta_x \colon \fX \to X$. The fact that the diameters of the clopen sets $\{U_{\ell}\}$ tend to zero, implies that $\Theta_x$ is a homeomorphism.  This is proved in detail in   \cite[Appendix~A]{DHL2016a}.
Moreover, $\Theta_x(x) =  e_{\infty} = (e\G_\ell) \in X$, the basepoint of the inverse limit \eqref{eq-invlimspace}, where $e\G_\ell = \G_\ell$ is the coset of the identity $e \in \G$.  
  Let $X$ have an ultrametric  metric such  that $\G$ acts on $X$ by isometries, for instance, let 
   \begin{equation}\label{dxmetric}d_X\left((x_\ell),(y_\ell) \right) = 2^{-m}, \quad \textrm{ where } m = \max \{\ell \mid x_\ell = y_\ell, \, \ell \geq 0\}. \end{equation}
  Then let $\dX$ be the ultrametric metric on $\fX$ induced from $d_X$ by the homeomorphism $\Theta_x$.
The minimal equicontinuous action $(X, \G, \Phi_x)$   is called the \emph{odometer model} centered at $x$ for  $(\fX,\G,\Phi)$. 

The group chain $\cG^x_\cU$ depends on $x$ and $\cU$, and one can introduce an equivalence relation which, for a given group $\G$, identifies the class of group chains with topologically conjugate associated odometer models. We   refer the interested reader to \cite{DHL2016a}.

\subsection{Unique ergodic invariant measure} \label{subsec-measure}

Given   a Cantor  action   $(\fX,\G,\Phi)$,  choose an adapted neighborhood basis $\cU$ and consider the corresponding group chain $\cG_\cU^x$ and the odometer model \eqref{eq-invlimspace}. The group $\G$ acts transitively on the coset space $X_{\ell} = \G/\G_{\ell}$ and we define a $\G$-invariant probability measure $\mu_{\ell}$ on $X_{\ell}$ by giving equal weight to each point (coset) in $X_{\ell}$.  Thus one has
  \begin{align}\label{eq-mubern}\mu_{\ell} (h \G_{\ell}) = \frac{1}{|\G: \G_{\ell}|}, \textrm{ for all }h\G_{\ell} \in X_{\ell} \textrm{ and all }\ell \geq 0,\end{align}
where $|\G: \G_{\ell}|$ denotes the index of $\G_{\ell}$ in $\G$. The unique $\G$-invariant measure on the inverse limit $X$ is defined as the limit of the pull-backs of these measures under the projection maps $X \to X_{\ell}$. Then the invariant measure $\mu$ on $\fX$ is the pull-back via the homeomorphism $\Theta_{x} \colon \fX \to X$.

Alternately, consider the closure $E = \overline{\Phi(\G)} \subset \Homeo(\fX)$ in the uniform topology. It  is a profinite compact group, called the \emph{Ellis} or \emph{enveloping group} \cite{Auslander1988,Ellis1969}. (If the action   $(\fX,\G,\Phi)$ is not assumed to be equicontinuous, then $\overline{\Phi(\G)}$ is only a semi-group.)  The group $E$ acts on $\fX$ with    isotropy group $E_x = \{g \in E \mid g\cdot x = x\}$ for $x \in \fX$. Then $E_x$ is a closed subgroup of $E$, and we have $\fX = E/E_x$. The profinite group $E$ has a unique Haar measure  $\widehat{\mu}$, which is invariant with respect to the action of $E$ on itself. The measure $\widehat{\mu}$ on $E$ pushes down to the measure $\mu$ on $\fX$.

\subsection{Lebesgue density theorem}

Let $(\fX,\G,\Phi,\mu)$ be a Cantor action   probability measure $\mu$ and ultrametric $d_\fX$ induced from \eqref{dxmetric}.  Denote by 
  $B(x,\epsilon) = \{ y \in \fX \mid d_\fX(x,y) < \epsilon\}$  the open ball with center $x$ of radius $\e >0$. The proof of the \emph{Lebesgue Density Theorem} in the formulation below can be found, for instance, in  \cite[Proposition 2.10]{Miller2008}.

\begin{theorem} \label{thm-polishlebesgue}
Let $\fX$ be a Polish space, and suppose $\fX$ has an ultrametric $d_\fX$ compatible with its topology. Let $\mu$ be a probability measure on $\fX$, and let $A$ be a Borel set of positive measure. Then the \emph{Lebesgue density} of $x$ in $A$, given by
  \begin{align}\label{eq-lebesgue}\lim_{\epsilon \to 0} \frac{\mu(A \cap B(x,\epsilon))}{\mu(B(x,\epsilon))}\end{align}
exists and is equal to $1$ for $\mu$-almost every $x \in A$.
\end{theorem}

We   give an important consequence of the Lebesgue Density Theorem.
 \begin{lemma}\label{lem-density}
 Let $(\fX, \G, \Phi,\mu)$  be a   Cantor action, with  invariant probability measure $\mu$.   
Assume there exists an element $g \in \G$ for which $\fX_g^{hol}$ has positive $\mu$-measure. 
Then there exists $x \in \fX_g^{hol}$ such that, for all $0 < \e <1$, there exists an adapted set $U_{\e}$ with $x \in U_{\e}$ and 
$\displaystyle \mu(U_{\e}  \cap \fX_g^{hol}) \geq (1-\e) \cdot  \mu(U_{\e})$.
\end{lemma}

\proof Since $\fX_g^{hol}$ has positive $\mu$-measure, by the Lebesgue Density Theorem \ref{thm-polishlebesgue} there exists a point $x \in \fX_g^{hol}$ of full Lebesgue density. 
For this point, choose an adapted neighborhood basis $\cU = \{U_{\ell} \subset \fX  \mid \ell \geq 0\}$ at  $x$. 
By the convergence of the limit in Theorem \ref{thm-polishlebesgue} there exists $\ell_{\e}$ so that 
 $\displaystyle \mu(U_{\ell} \cap \fX_g^{hol}) \geq (1-\e) \cdot  \mu(U_{\ell})$ for $\ell \geq \ell_0$. Then set $U_{\e} = U_{\ell_{\e}}$.
\endproof

\subsection{Locally quasi-analytic}\label{subsec-lqa}

 The quasi-analytic property for Cantor actions  was introduced by
 {\'A}lvarez L{\'o}pez and  Candel  in  \cite[Definition~9.4]{ALC2009} as a generalization of the notion of a \emph{quasi-analytic action} studied by Haefliger for actions of pseudogroups of real-analytic diffeomorphisms in \cite{Haefliger1985}.  The authors introduced a local form of the quasi-analytic property  in \cite{HL2018a,HL2018b}: 
  
\begin{defn} \cite[Definition~2.1]{HL2018b} \label{def-LQA} A topological action       $(\fX,\G,\Phi)$ on a metric Cantor space $\fX$,  is   \emph{locally quasi-analytic}  if there exists $\e > 0$ such that for any non-empty open set $U \subset \fX$ with $\diam (U) < \e$,  and  for any non-empty open subset $V \subset U $, and elements $g_1 , g_2 \in \G$
 \begin{equation}\label{eq-lqa}
  \text{if the restrictions} ~~ \Phi(g_1)|V = \Phi(g_2)|V, ~ \text{ then}~~ \Phi(g_1)|U = \Phi(g_2)|U. 
\end{equation}
 The action is said to be \emph{quasi-analytic} if \eqref{eq-lqa} holds for $U=\fX$.
\end{defn}

In other words, $(\fX,\G,\Phi)$ is locally quasi-analytic if for every $g \in \G$, the homeomorphism $\Phi(g)$ has unique extensions on the sets of diameter $\e>0$ in $\fX$, with $\e$ uniform over $\fX$. We note that an effective action $(\fX,\G,\Phi)$ is topologically free if and only if it is quasi-analytic \cite[Proposition 2.2]{HL2018b}. 

Recall that a group $\G$ is \emph{Noetherian}  \cite{Baer1956} if every increasing chain of subgroups has a maximal element. Equivalently, a group is Noetherian if every subgroup of $\G$ is finitely generated.

\begin{thm}\label{thm-noetherian} \cite[Theorem~1.6]{HL2018b}
Let    $\G$ be a Noetherian group. Then   a  minimal equicontinuous Cantor action $(\fX,\G,\Phi)$   is locally quasi-analytic.
\end{thm}
 A finitely-generated nilpotent group is Noetherian, so as a corollary we obtain that all Cantor actions by finitely-generated nilpotent groups are locally quasi-analytic.   
  Examples of locally quasi-analytic actions which are not quasi-analytic are easy to construct; see for instance \cite[Example A.4]{HL2018b}.

\section{Invariance}\label{sec-invariance}

 We recall   notions of equivalence for Cantor actions, considered as topological dynamical systems. For each notion considered, we show that the measurable dynamical systems property that an action has essential holonomy is preserved for equivalent actions. The work \cite{GPS2019} gives an excellent comparison of the various notions of equivalence for the case of Cantor actions by $\G = \mZ^n$.
 
 First, we recall the most basic equivalence of actions.

  \begin{defn} \label{def-isomorphism}
  Cantor actions $(\fX_1, \G_1, \Phi_1)$ and $(\fX_2, \G_2, \Phi_2)$   are said to be \emph{isomorphic}, or \emph{conjugate},  if there is a homeomorphism $h \colon \fX_1 \to \fX_2$ and a group isomorphism $\Theta \colon \G_1 \to \G_2$ so that 
\begin{equation}\label{eq-isomorphism}
\Phi_1(\gamma) = h^{-1} \circ \Phi_2(\Theta(\gamma)) \circ h   \in   \Homeo(\fX_2) \  {\rm for \ all} \ \gamma \in \G_1 \ .
\end{equation}
 \end{defn}

\begin{prop}\label{prop-isoessholo}
Let  $(\fX_1, \G_1, \Phi_1)$ and $(\fX_2, \G_2, \Phi_2)$ be isomorphic Cantor actions. Then $(\fX_1, \G_1, \Phi_1)$ has essential holonomy if and only if  $(\fX_2, \G_2, \Phi_2)$ has essential holonomy.
\end{prop}
\proof
Let  $(\fX_1, \G_1, \Phi_1)$ and $(\fX_2, \G_2, \Phi_2)$ be isomorphic Cantor actions, with   map  $h \colon \fX_1 \to \fX_2$ and isomorphism $\Theta \colon \G_1 \to \G_2$. Let $\mu_1$ be the unique invariant probability measure for $(\fX_1, \G_1, \Phi_1)$, and $\mu_2$  be the unique invariant probability measure for $(\fX_2, \G_2, \Phi_2)$. Then $h^*(\mu_2)$  is an invariant probability measure for $(\fX_1, \G_1, \Phi_1)$, hence $\mu_1 = h^*(\mu_2)$. 
Moreover, for $g \in \Gamma_1$ the action $\Phi_1(g)$ has essential holonomy if and only if $\Phi_2(\Theta(g))$ has essential holonomy, so 
$h(\fX_{\G_1}^{hol}) =  \fX_{\G_2}^{hol}$, and the claim follows.
\endproof

The notion of \emph{return equivalence}  is  the analog for Cantor actions of  Morita equivalence for $C^*$-algebras. This equivalence is  weaker than the notion of isomorphism, and is natural when considering the Cantor systems defined by the holonomy actions for matchbox manifolds, as   in    \cite{HL2018a,HL2018b}.

 \begin{defn}\label{def-re}
 Cantor actions $(\fX_1, \G_1, \Phi_1)$ and $(\fX_2, \G_2, \Phi_2)$    are said to be  \emph{return equivalent} if there exist non-empty clopen subsets $U_i \subset \fX_i$, for $i=1,2$,  such that $U_i$ is adapted to the action $\Phi_i$, and there is a homeomorphism $h \colon U_1 \to U_2$ whose induced homomorphism $h_* \colon \Homeo(U_1) \to \Homeo(U_2)$  restricts to an isomorphism  $\Theta \colon \cH_{U_1} \to \cH_{U_2}$.
 \end{defn}
Note that when $U_i = \fX_i$ and both actions are effective, then this definition reduces to the usual notion of isomorphism of the actions, with induced group isomorphism $\Theta \colon \G_1 \cong \cH_{\fX_1} \to \cH_{\fX_2} \cong \G_2$.  

\begin{prop}\label{prop-REessholo}
Let  $(\fX_1, \G_1, \Phi_1)$ and $(\fX_2, \G_2, \Phi_2)$ be return equivalent Cantor actions. Then $(\fX_1, \G_1, \Phi_1)$ has essential holonomy if and only if  $(\fX_2, \G_2, \Phi_2)$ has essential holonomy.
\end{prop}
\proof
Assume that  $(\fX_1, \G_1, \Phi_1)$ has essential holonomy, and let $U_1 \subset \fX_1$ and $U_2 \subset \fX_2$ be clopen sets such that the restricted actions are isomorphic by $\Theta \colon   \cH_{U_1} \to \cH_{U_2}$. Let $g \in \G_1$ be such that 
$\fX_g^{hol}$ has positive $\mu_1$-measure, and thus there exists $x \in \fX_1$ such that $x$ is fixed by $\Phi_1(g)$ and $ \fX_{1,g}^{hol}$ has Lebesgue density 1 at $x$. The action $\Phi_1$ is minimal on $\fX_1$ so there exists $k \in \G_1$ such that $k \cdot x \in U_1$. Then $g' = kgk^{-1} \in \G_1$ has fixed point $x' = kx$ which is a point of Lebesgue density 1 in the set 
$\fX_{1,g'}^{hol}$.  
  As $x' \in U_1 \cap \Phi_1(g')(U_1)$ and $U_1$ is adapted, we have $U_1 = \Phi_1(g')(U_1)$ and so $g' \in \G_{1,U}$. 
  
  The action of $\G_{1,U_1}$ on $U_1$ is minimal, so the renormalized    measure $\mu_1' = \mu(U_1)^{-1} \mu_1 \mid U_1$ is the unique invariant probability measure for the restricted action of $\G_{1,U_1}$ on $U_1$.
  
The set $\fX_{1,g'}^{hol} \cap U_1$ has Lebesgue density 1 at $x'$,  hence its image $h(x') \in h(\fX_{1,g'}^{hol} \cap U_1) \subset U_2$   is also a point of Lebesgue density 1 for the   action of $\Theta(g')$ on $U_2$, with corresponding renormalized measure $\mu_2'$ on $U_2$. Thus,  $(\fX_2, \G_2, \Phi_2)$ has essential holonomy. The converse follows similarly. 
\endproof

The  notion of \emph{continuous orbit equivalence} for   Cantor actions  was introduced in     \cite{BoyleTomiyama1998}. It  is the analogue for topological dynamics of measurable orbit equivalence for measurable actions, as first introduced by Dye \cite{Dye1959}.   Continuous orbit equivalence plays a fundamental role in   the classification of group actions on Cantor sets (see for example \cite{Renault2008}).

 \begin{defn}\label{def-coe}
 Cantor actions  $(\fX_1, \G_1, \Phi_1)$ and $(\fX_2, \G_2, \Phi_2)$  are said to be \emph{continuously orbit equivalent} if there exists a homeomorphism $h \colon \fX_1 \to \fX_2$ and continuous functions
\begin{enumerate}
\item $\alpha \colon G_1 \times \fX_1 \to G_2$,  $h(\Phi_1(g_1)(x_1)) = \Phi_2(\alpha(g_1 , x_1))(h(x_1))$ for all   $g_1 \in G_1$ and $x_1 \in \fX_1$; 
\item $\beta \colon  G_2 \times \fX_2 \to G_1$,  $h^{-1}(\Phi_2(g_2)(x_2)) = \Phi_1(\beta(g_2, x_2))(h^{-1}(x_2))$ for all $g_2 \in G_2$ and $x_2 \in \fX_2$. 
\end{enumerate}
\end{defn}
 The homeomorphism  $h$ is called a  \emph{continuous orbit equivalence} between the two actions.
Note that  the functions $\alpha$ and $\beta$ are not assumed to satisfy the cocycle property.

We have the following result of Cortez and Medynets:
\begin{thm} \cite{CortezMedynets2016} \label{thm-CM}
Let  $(\fX_1, \G_1, \Phi_1)$ and $(\fX_2, \G_2, \Phi_2)$ be topologically free Cantor actions. If the actions $\Phi_1$ and $\Phi_2$ are continuously orbit equivalent,  then they  are return equivalent.  
\end{thm}

This result was generalized by the authors:

\begin{thm}\cite{HL2018b} \label{thm-main1}
Let  $(\fX_1, \G_1, \Phi_1)$ and $(\fX_2, \G_2, \Phi_2)$ be locally quasi-analytic Cantor actions. If the actions $\Phi_1$ and $\Phi_2$ are continuously orbit equivalent,  then they  are return equivalent.  
\end{thm}

For the general case (not necessarily locally quasi-analytic)  of Cantor actions, we have the following:
\begin{prop}\label{prop-COEessholo}
Let  $(\fX_1, \G_1, \Phi_1)$ and $(\fX_2, \G_2, \Phi_2)$ be  continuously orbit equivalent Cantor actions. Then $(\fX_1, \G_1, \Phi_1)$ has essential holonomy if and only if  $(\fX_2, \G_2, \Phi_2)$ has essential holonomy.
\end{prop}
\proof
 Let   $h \colon \fX_1 \to \fX_2$ be the homeomorphism given by  Definition~\ref{def-coe}. Let $\mu_2$ be the unique invariant probability measure for $(\fX_2, \G_2, \Phi_2)$, then the pull-back $\widetilde{\mu}_2 = h^*(\mu_2)$ is a probability measure on $\fX_1$. Condition (1) implies   $\widetilde{\mu}_2$ is invariant under the action $(\fX_1, \G_1, \Phi_1)$, hence $\mu_1 = \widetilde{\mu}_2$.

Assume that  $(\fX_1, \G_1, \Phi_1)$ has essential holonomy, and let $g \in \G_1$ be such that 
$\fX_g^{hol}$ has positive $\mu_1$-measure, so there exists $x \in \fX_1$ such that   $\Phi_1(g)(x) = x$ and $ \fX_g^{hol}$ has $\mu_1$-Lebesgue density 1 at $x$.

Let  $\alpha \colon G_1 \times \fX_1 \to G_2$ the map given by condition (1) so that   $h(\Phi_1(g_1)(x_1)) = \Phi_2(\alpha(g_1 , x_1))(h(x_1))$ for all   $g_1 \in G_1$ and $x_1 \in \fX_1$. 
Then for $g_1 =g$ there exists a clopen set $U_1 \subset \fX_1$ so that $x \in U_1$ and $\alpha(g, y) \in \G_2$ is continuous for $y \in U_1$. Since $\G_2$ is a discrete space, we can assume that $U_1$ is sufficiently small so that $g_2 = \alpha(g, y) \in \G_2$ is constant for   $y \in U$.

We then have $h(\Phi_1(g)(y)) = \Phi_2(g_2)(h(y))$ for all $y \in U_1$. Set $U_2 = h(U_1)$ then this states that $h_{U_1} \colon U_1 \to U_2$ conjugates the action of $\Phi_1(g_1)$ on $U_1$ with the action of $\Phi_2(g_2)$ on $U_2$. Thus, $\Phi_2(g_2)$ has non-trivial holonomy at the fixed point $z=h(x)$ and this is a point of Lebesgue density 1 for the action of  $\Phi_2(g_2)$. In particular, the action $(\fX_2, \G_2, \Phi_2)$ has essential holonomy. 

The converse implication is proved similarly, using condition (2) of Definition~\ref{def-coe}.
\endproof

\section{Dynamics and the lower central series}\label{sec-proof}
 
  Theorem \ref{thm-main0}   relates the non-trivial essential holonomy property for a     Cantor action $(\fX, \G, \Phi)$, with the  lower central series of $\G$. In this section, we show that if a Cantor action of $\G$ is locally quasi-analytic  and has essential holonomy, then every  commutator subgroup in the lower central series of $\G$ has elements with positive measure sets of points with non-trivial holonomy.

First, recall the construction of the lower central series for   $\G$.  
Set $\gamma_1(\G) = \G$, and  for  $i \geq 1$,    let \break $\gamma_{i+1}(\G) = [\G,\gamma_{i}(\G)]$ be the commutator subgroup, which  is a normal subgroup of $\G$. Then for $a \in \gamma_{i}(\G)$ and $b \in \gamma_{j}(\G)$   the commutator $[a,b] \in \gamma_{i+j}(\G)$. Moreover, these subgroups form a descending chain
 \begin{equation}\label{eq-filtration}
\G = \gamma_{1}(\G) \supset \gamma_{2}(\G) \supset \cdots  \supset  \gamma_{n}(\G) \supset \cdots  \ .
\end{equation}
  Note that each quotient group $\gamma_{i}(\G)/\gamma_{i+1}(\G)$ is abelian. 
  
  The group $\G$ is nilpotent of length $n_{\G}$ if there exists an $n_\G > 0$ such that $\gamma_{n}(\G)$ is the trivial group for $n > n_{\G}$, and $\gamma_{n_{\G}}(\G)$ is non-trivial. It follows that  every  element in $\gamma_{n_{\G}}(\G)$ commutes with every element of $\G$; that is, $\gamma_{n_{\G}}(\G)$ is contained in the center of $\G$. 

Denote by $\Phi_n \colon \gamma_{n}(\G) \times \fX \to \fX$ the restriction of the action $\Phi$ to the subgroup $\gamma_{n}(\G)$.
  
\begin{defn}\label{def-turbdepth}
A   Cantor action $(\fX, \G, \Phi)$, with invariant probability measure $\mu$, is said to have \emph{essential holonomy at depth $n_t$} if the restricted action   $(\fX, \gamma_{n_t}(\G), \Phi_{n_t})$ has essential holonomy, but $(\fX, \gamma_{n}(\G), \Phi_n)$ has no essential holonomy for $n > n_t$. The action has \emph{essential holonomy at infinite depth} if for all $n \geq 1$, the restricted action $(\fX, \gamma_n(\G), \Phi_n)$ has essential holonomy.
\end{defn}

  Here is our main technical result, which combined with Theorem \ref{thm-noetherian} implies Theorem~\ref{thm-main0}.
  \begin{prop}\label{prop-commutators}
Let $(\fX, \G, \Phi)$ be a   Cantor action which is locally quasi-analytic. If the action   has essential holonomy, then it has essential holonomy of infinite depth.
\end{prop}
\proof
Suppose that $(\fX, \G, \Phi)$  has essential holonomy at depth $n_t$.  We show that this leads to a contradiction by localizing the action to a sufficiently small adapted set. Recall that $\mu$ denotes the unique invariant probability measure for the action.

 The set 
$\fX_{\gamma_{(n_t)}(\G)}^{hol}$ has positive $\mu$-measure, while the set $\fX_{\gamma_{(n_t +1)}(\G)}^{hol}$ has   $\mu$-measure zero.
Thus, there  exists $g \in \gamma_{n_t}(\G)$ such that $\mu(\fX_{g}^{hol}) > 0$. Moreover, for all $\tau \in \gamma_{(n_t+1)}(\G)$ we have $\mu(\fX_{\tau}^{hol}) = 0$.

      Let $x \in \fX$ be such that $g \cdot x =x$, and $\fX_g^{hol}$ has   Lebesgue density $1$ at $x$. Let $U \subset \fX$ be an adapted clopen set with $x \in U$, and so $g \in \G_U$,  and sufficiently small diameter such that the restricted action $\Phi_U \colon \G_U \times U \to U$ is quasi-analytic, and we have $\mu(\fX_g^{hol} \cap U) \geq 3/4 \cdot \mu(U)$, as in  Lemma \ref{lem-density}.

  By the assumption that $x \in \fX_g^{hol}$, there exist  $y \in U$ so that $z = g\cdot y \ne y$. Note that   $z \ne x$ as $g \cdot x = x$. Let $V \subset U$ with $y \in V$ be a sufficiently small clopen set such that $(g \cdot V) \cap V = \emptyset$. Then choose $k_y \in \G_U$ such that $k_y \cdot x \in V$. Then replace $y$ with $k_y \cdot x$ and we have $g \cdot y \ne y$.

Let $\tau = [g, k_y]$ be the commutator. Then $g \in  \gamma_{n_t}(\G)$ implies that $\tau \in \gamma_{(n_t+1)}(\G) \cap \G_U$.
By definition we have $g \cdot k_y = \tau \cdot k_y \cdot g$.
Then for $w \in \fX_g^{hol} \cap U$ calculate
\begin{equation}\label{eq-gamma}
g \cdot (k_y \cdot w) =  \tau \cdot k_y \cdot g \cdot w =  \tau \cdot (k_y \cdot  w) \ .
\end{equation}
We use the identity \eqref{eq-gamma} to prove the key observation:
\begin{lemma}\label{lem-disjoint}
The sets $k_y\cdot (\fX_g^{hol} \cap U)$ and  $(\fX_g^{hol} \cap U)$ are $\mu$-a.e. disjoint.  
\end{lemma}
\proof
We must show that $\mu\left( k_y\cdot (\fX_g^{hol} \cap U) \cap (\fX_g^{hol} \cap U)\right) = 0$.

Suppose that  $w \in \fX_g^{hol} \cap U$  satisfies $k_y\cdot w \in \fX_g^{hol} \cap U$. Then $k_y \cdot w$ is a fixed-point for the action of $g$, and so by \eqref{eq-gamma} we have 
$k_y \cdot w $ is a fixed point for the action of $\tau = [g, k_y]$.

 As  $g \in \gamma_{n_t}(\G)$, we have that $\tau \in \gamma_{(n_t+1)}(\G)$. Then by the definition of the index $n_t$, for $\mu$-almost all  $w \in \fX_g^{hol} \cap U$, we have that $\Phi_U(\tau)$ has trivial germinal holonomy at $k_y \cdot w$, as the action of $\Phi(k_y)$ is a measure preserving homeomorphism.

 Thus, there exists a   clopen set $W_{w}  \subset U$ with  $k_y \cdot w \in W_{w}$ such that $\Phi_U(\tau)$ acts as the identity on $W_{w}$. As we chose $U$ so that the action $\Phi_U$ is quasi-analytic on $U$, this implies that   $\Phi_U(\tau)$ acts as the identity on $U$.
 However, we also have $g \cdot y \ne y$, so using the identity  \eqref{eq-gamma} again, we have
 $\tau \cdot y = \tau \cdot (k_y \cdot  x) \ne k_y \cdot x$, and thus $\Phi_U(\tau)$ does not act as the identity on $U$, which is a contradiction.
It follows that for   $\mu$-almost all  $w \in \fX_g^{hol} \cap U$ we have $k_y \cdot w \not\in \fX_g^{hol} \cap U$. Hence  $\mu\left( k_y\cdot (\fX_g^{hol} \cap U) \cap (\fX_g^{hol} \cap U)\right) = 0$.
 \endproof

We now complete the proof of  Proposition~\ref{prop-commutators}.  As $\mu$ is invariant under the action of $\Phi_U$ we have  $\mu\left(k_y \cdot (\fX_g^{hol} \cap U) \right) = \mu\left(\fX_g^{hol} \cap U\right)\geq 3/4 \cdot \mu(U)$. But then by Lemma~\ref{lem-disjoint} we obtain the contradiction 
\begin{equation}\label{eq-toomuch}
\mu(U) \geq  \mu\left( k_y\cdot (\fX_g^{hol} \cap U)\right) + \mu \left(\fX_g^{hol} \cap U\right) \geq (3/4 + 3/4) \mu(U) > \mu(U) \ .
\end{equation}
 Thus, the action $(\fX, \G, \Phi)$ cannot have essential holonomy at finite depth.
 \endproof

 \begin{cor}\label{cor-commutators}
Let $(\fX, \G, \Phi)$ be a  Cantor action, and suppose that $\G$ is a finitely-generated nilpotent group. Then the  set of points with non-trivial holonomy has measure 0;  that is, $(\fX, \G, \Phi)$ does not have essential holonomy. 
\end{cor}
 \proof
 By Theorem \ref{thm-noetherian}, the action $(\fX, \G, \Phi)$ is locally quasi-analytic. As  $\G$ is nilpotent,   the commutator group $\gamma_{n_{\G}}(\G)$ is in the center of $\G$,   hence the  set of points with non-trivial holonomy for the action of $\gamma_{n_{\G}}(\G)$ on $\fX$ has measure 0;  that is, it has no essential holonomy.  Thus by Proposition~\ref{prop-commutators} the action of $\G$ on $\fX$ has no essential holonomy.
 \endproof
 
 \section{Commutators and not locally quasi-analytic actions}
 
 In this section, we exhibit  a family of examples to show that the assumption that an action $(\fX,\G,\Phi)$ is locally quasi-analytic is essential for the conclusions of Theorem \ref{thm-main0} and Corollary \ref{cor-commutator}.
 
 \begin{thm}\label{thm-examples}
 Let $\{n_\ell\}_{\ell \geq 1}$, $n_\ell \geq 6$ be a sequence of positive numbers, such that $n_\ell > 2 n_{\ell-1}$. Let $S_{\ell}$ be a set with $n_\ell$ elements, and let $\fX = \prod_{\ell \geq 1} S_\ell$. Let $A_\ell$ be the alternating group on $n_\ell$ symbols. There exists a countably generated subgroup $\G \subset \prod_{\ell \geq 1} A_\ell$ with the following properties.
 \begin{enumerate}
 \item The lower central series of $\G$ stabilizes, i.e. $\gamma_n(\G) = [\G,\G]$ for $n \geq 2$.
 \item The minimal equicontinuous action $(\fX,\G,\Phi)$ is not locally quasi-analytic and has essential holonomy, i.e. the set of points with non-trivial holonomy has full measure.
 \item The induced action $(\fX,[\G,\G],\Phi)$ of the commutator subgroup is minimal, not locally quasi-analytic and it has no essential holonomy. 
 \end{enumerate}
\end{thm}

At the moment we are not aware of an action of a finitely generated group which exhibits similar properties for the commutator action. 

The family of examples in Theorem \ref{thm-examples} is an amalgam of the family considered in \cite[Theorem 1.5]{ABLN2022} and \cite[Theorem 6.1]{GL2019}. The idea to use actions of alternating groups comes from \cite[Theorem 1.5]{ABLN2022}, and the construction of an element with positive measure set of points with holonomy is the same as in \cite[Theorem 6.1]{GL2019}. When constructing the commutator subgroup, we have to restrict to using the direct sum of alternating groups instead of the direct product to ensure that the action of the commutator in our family has no essential holonomy (in fact, no points with non-trivial holonomy at all), and consequently the group we construct is not finitely generated.
Merging these constructions   yields the actions in Theorem~\ref{thm-examples}.

\proof
The demonstration of Theorem~\ref{thm-examples} follows from Lemmas~\ref{lem-construction1}, \ref{lem-construction2} and \ref{lem-construction3}  below.
\begin{lemma}\label{lem-construction1}
Consider the direct sum 
\begin{equation}\label{eq-defG}
G =  \bigoplus_{\ell \geq 1} A_\ell = \{g = (g_1,g_2,\ldots ) \mid g_\ell \in A_\ell, \quad g_\ell = e_\ell \textrm{ for all but a finite set of }\ell \} ~ ,
\end{equation}

 where $e_\ell$ is the identity element in $A_\ell$. Then the component-wise action of $G$ on $\fX$ is minimal, equicontinuous, not locally quasi-analytic and has no essential holonomy.
 \end{lemma}
 
 \proof 
 The group $G$ acts on the direct product space $\fX$ minimally, since the action of every $A_\ell$ on $S_\ell$ is transitive. Define the metric on $\fX$ by setting $d_{\fX}(x,y) = 1$ for all $x=(x_i),y=(y_i) \in \fX$   if $x_1 \ne y_1$, and otherwise
   $$d_{\fX}(x,y) = 2^{-K}, \quad K= \max \{k \geq 1 \mid x_\ell = y_\ell \textrm{ for all } 1 \leq \ell \leq k\}.$$
 Since $G$ acts on each component of $\fX$ by permutations, its action is equicontinuous.
 
 If $g \in G$ has a fixed point, then this point has a clopen neighborhood entirely fixed by the action of $g$, since $g$ acts non-trivially on at most a finite number of factors in the direct product $\fX$. Thus the action of $G$ on $\fX$ has no essential holonomy. 
 
 We show that the action of $G$ is not locally quasi-analytic. Label the elements in $S_4$ by $\{1,2,3, \ldots, n_\ell\}$, and let $a_\ell = (123) \in A_\ell$, i.e. $a_\ell$ fixes all symbols except $1$, $2$ and $3$. Let $g \in G$, then there is $m \geq 1$ such that $g_\ell = a_\ell$ for a subset of $1 \leq \ell \leq m$, and otherwise $g_\ell = e_\ell$. Then $g$ fixes every point in the clopen set $V = \prod_{1 \leq \ell \leq m} \{4,\ldots,n_\ell\} \times \prod_{\ell > m} S_\ell$. It is clear that the restriction $g|V$ has multiple extensions to larger sets, for instance, take $\widetilde g$ which acts non-trivially on the same factors as $g$ by the permutation $b_\ell = (132)$.
 \endproof
 
For sequences $\{n_\ell\}_{\ell \geq 1}$ with the additional property that $n_{\ell+1} > 2 n_{\ell}$ for $\ell \geq 1$ we now realize $G$ as the commutator of a discrete group $\G$ whose action has essential holonomy. 

For $\ell \geq 1$, define a permutation $\gamma_\ell = (n_\ell-1\, n_\ell)$, i.e. $\gamma_\ell$ fixes every vertex in $S_\ell$ except the last two, and let $\gamma = (\gamma_1, \gamma_2,\ldots) \in \prod_{\ell \geq 1} A_\ell$ be an element in the \emph{direct product} of alternating groups. Define
 \begin{equation}\label{eq-defGamma}
\G = \langle \gamma, \G \rangle \subset \prod_{\ell \geq 1} A_\ell \ . 
\end{equation}

 \begin{lemma}\label{lem-construction2}
The commutators $  [G,G] = [\Gamma, G] = [\Gamma, \Gamma] = G$. 
\end{lemma}
\proof
 
For each $\ell \geq 1$, we have    $n_{\ell} \geq 5$, hence the alternating group $A_{\ell}$ is simple and thus perfect, and  thus  $G$ is also perfect; that is, $G = [G,G]$.

Next we show that $[\Gamma, \Gamma] = G$. Note that the restriction   $\gamma|S_\ell = \gamma_\ell$ is an odd permutation, and $\gamma$ acts on the components of $\prod_{\ell \geq 1}A_\ell$ independently, so for any $g \in G$ the commutator  $[\gamma, g]|S_\ell$ is an even permutation, hence $[\gamma, g] \in G$. Thus, $G = [G,G] \subset [\G, G] \subset [\G,\G] \subset G$. \endproof

Using the metric $d_\fX$, defined above, it is convenient to think of the finite products $\prod_{1 \leq \ell \leq k}S_\ell$ as sets of vertices of a rooted tree at level $k \geq 1$, where the root is a single vertex at level $0$ (it is omitted from $\fX$). In such a tree, each vertex at level $\ell$ is connected to $n_{\ell+1} = |S_{\ell+1}|$ vertices at level $\ell+1$, and so one can think of vertices in $\prod_{1 \leq \ell \leq k}S_\ell$ as labeled by finite words $s_1 \cdots s_k$, where $s_\ell \in S_\ell$. An element of the direct product $\gamma \in \prod_{\ell \geq 1} A_\ell$ acts on the tree so that, for a given $s_\ell \in S_\ell$,  $\gamma \cdot s_\ell = \gamma_\ell \cdot s_\ell$, and this action depends only on the symbol $s_\ell$ in the finite word $s_1 \cdots s_k$. An element of $\fX$ is a infinite sequence $s_1 s_2 \cdots$, where each truncated word $s_1 \cdots s_k$ corresponds to a vertex in $\prod_{1 \leq \ell \leq k}S_\ell$.

Open balls of diameter $2^{-K}$ in the metric $d_\fX$ are the sets of infinite words in $\fX$ which coincide for their first $K$ symbols. The measure of each such open ball is $$\frac{1}{|\prod_{1 \leq \ell \leq K}S_\ell|} = \frac{1}{n_1 \cdots n_K}.$$

We now show that $\gamma$ has positive measure set of points with non-trivial holonomy by explicitly computing the Lebesgue density of this set at each point. The argument is the same as in \cite[Theorem 6.1]{GL2019} and we give it here for completeness and for the convenience of the reader.

\begin{lemma}\label{lem-construction3}
Let ${\rm Fix}(\gamma)$ be the set of fixed points of the element $\gamma \in \prod_{\ell \geq 1} A_\ell$ defined above. Then every point in ${\rm Fix}(\gamma)$ has non-trivial holonomy, and $\mu({\rm Fix}(\gamma)) > 0$.
\end{lemma}

\proof
First note that $\gamma$ fixes an infinite path $s = s_1 s_2 \cdots \in \fX$ if and only if for all $\ell \geq 1$ we have $s_\ell \ne n_\ell$ and $s_\ell \ne n_\ell - 1$. We claim that each such point has non-trivial holonomy. Indeed, let $V$ be an open neighborhood of $s$. Then there is an $m_V \geq 1$ such that 
  $$V \cap \left(\prod_{1 \leq \ell \leq m_V} \{s_\ell\} \times \prod_{ \ell > m_V} S_\ell \right) =\prod_{1 \leq \ell \leq m_V}\{ s_\ell \} \times \prod_{ \ell > m_V} S_\ell,$$
 and so $V$ contains the points which are moved by the action of $\gamma$. 
 
 For each clopen ball $U_\ell$ around a fixed point $s$, the action of $\gamma$ permutes two clopen balls in $U_\ell$ consisting of sequences starting with $s_1 \cdots s_\ell (n_{\ell+1}-1)$ and $s_1 \cdots s_\ell n_{\ell+1}$, which have the total measure $2/(n_1 \cdots n_{\ell+1})$. Each of the remaining $n_{\ell+1} - 2$ clopen balls determined by words of length $\ell+1$ contains $2$ subsets of sequences starting with $s_1 \cdots s_{\ell+1} (n_{\ell+2}-1)$ and $s_1 \cdots s_{\ell+1} n_{\ell+2}$ permuted by the action, whose total measure is $2(n_{\ell+1} - 2)/(n_1 \cdots n_{\ell+1} n_{\ell+2})$. Continuing by induction, we compute the upper bound on the measure of the complement of the set ${\rm Fix}(\gamma)$:
 $$\mu(U_\ell - {\rm Fix}(\gamma)) = \frac{1}{n_1 \cdots n_\ell} \left( \frac{2}{n_{\ell+1}} + \frac{2(n_{\ell+1} - 2)}{n_{\ell+1} n_{\ell+2}} +  \frac{2(n_{\ell+1}-2)( n_{\ell+2} - 2)}{n_{\ell+1} n_{\ell+2} n_{\ell+2}} + \cdots \right) <\frac{1}{n_1 \cdots n_\ell} \sum_{i \geq 1} \frac{1}{n_{\ell+i}}. $$
 Since we assume that $n_{\ell+i} > 2 n_{\ell+i-1} > 2^{i-1} n_{\ell+1}$, we obtain that
   $$\mu(U_\ell - {\rm Fix}(\gamma)) < \frac{1}{n_1 \cdots n_\ell} \frac{4}{n_{\ell+1}},$$
 and so 
   $$\mu(U_\ell \cap {\rm Fix}(\gamma)) > \frac{1}{n_1\cdots n_\ell} - \frac{4}{n_1 \cdots n_{\ell+1}}.$$  
   It follows that for every point in ${\rm Fix}(\gamma)$ the Lebesgue density is $1$, namely
  \begin{equation}
1 =  \lim_{\ell \to \infty} \left(1 - \frac{4}{n_{\ell+1}} \right) \leq \lim_{\ell \to \infty} \frac{\mu (U_\ell \cap {\rm Fix}(\gamma))}{\mu(U_\ell)} \leq  1 \ .
\end{equation}
    Thus $\gamma$ has the positive measure set of points with non-trivial holonomy. \endproof
      
    We have shown that  the action of $\G$ on $\fX$ has essential holonomy, while the action of its commutator $[\G,\G] = G$ has no essential holonomy, which proves the assertions of Theorem~\ref{thm-examples}.
\endproof
 

\end{document}